\documentclass[reqno]{amsart}
\usepackage{hyperref}
\newtheorem{theorem}{Theorem}[section]
\newtheorem{definition}{Definition}[section]

\newtheorem{lemma}{Lemma}[section]

\newtheorem{remark}{Remark}[section]

\AtBeginDocument{{\noindent\small
\emph{Journal of Fractional Calculus and Applications}\\
Vol. 9(1) Jan.   2018,  pp. 120-136.\newline
ISSN: 2090-5858.\\
http://fcag-egypt.com/Journals/JFCA/\\
------------------------------------------------------------------------------------------------}
\vspace{9mm}}

\begin{document}
\title[\hfilneg JFCA-2018/9(1)\hfil Existence of solutions for boundary conditions]
{ Existence of solutions to fractional differential equations with
three-point boundary conditions at resonance with general conditions}

\author[Bin-Bin He\hfil JFCA-2018/9(1)\hfilneg]
{Bin-Bin He}

\address{Bin-Bin He \newline
College of Information Science and Technology,  Donghua University,  Shanghai,  P.R. China}
\email{hebinbin45@126.com}

\thanks{Submitted Mar. 19, 2017.  Revised May 13, 2017. }

\setcounter{page}{120}

\maketitle

\begin{abstract}
This paper presents a new technique to investigate the existence of solutions
to  fractional three-point boundary value problems at resonance
in a Hilbert space.  Based on the proposed method, the
restricted conditions $A^2\xi^{2\alpha-2}=A\xi^{\alpha-1}$ and  $A^2\xi^{2\alpha-2}=I$
on the operator $A$, which have been used in \cite{Zhou2}, are removed.
It is shown that the system under consideration admits at least one solution by applying
coincidence degree theory. Finally, an illustrative example is presented.
\end{abstract}

\section{Introduction}
In this article, we consider the problem of the existence of solutions for
the following fractional three-point boundary value problems(BVPs)
at resonance\\
\begin{equation}\label{model-equ}
\begin{gathered}
D_{0^+}^{\alpha}x(t)=f(t,x(t),D_{0^+}^{\alpha-1}x(t)), \quad 1<\alpha\leq 2,\;
  t\in [0,1], \\
I_{0+}^{2-\alpha}x(t)|_{t=0}=\theta, \quad x(1)=Ax(\xi),
\end{gathered}
\end{equation}
where  $D_{0^+}^{\alpha}$ and  $I_{0^+}^{\alpha}$ represents the Riemann-Liouville
differentiation and integration, respectively;
 $\theta$ is the zero vector in
 $l^2:=\{x=(x_1,x_2,\dots,.)^{\top}:\sum_{i=1}^{\infty}|x_i|^2<\infty\}$;
$A:l^2\to l^2$ is a bounded linear operator satisfying
$1\leq\operatorname{dim}\ker (I-A\xi^{\alpha-1})<\infty$; $\xi\in(0,1)$
 is a fixed constant; $f:[0,1]\times l^2\times l^2\to l^2$ is a Carath\'eodory
function, that is,
\begin{itemize}
\item[(i)] for each $(u,v)\in l^2\times l^2$, $t\mapsto f(t,u,v)$ is measurable
on $[0,1]$;

\item[(ii)] for a.e. $t\in[0,1],\;(u,v)\mapsto f(t,u,v)$ is continuous on
$l^2\times l^2$;

\item[(iii)] for every bounded set $\Omega\subseteq  l^2\times l^2$,
$\{f(t,u,v):(u,v)\in \Omega\}$ is a relatively compact  set in $l^2$.
Moreover,
\[
\varphi_\Omega(t)=\sup\{\|f(t,u,v)\|_{l^2}:(u,v)\in \Omega\}
\in L^1[0,1],
\]
where  $\|x\|_{l^2}=\sqrt{\sum_{i=1}^{\infty}|x_i|^2}$ is the norm
of $x=(x_1,x_2,\dots,\cdot)^\top$ in $l^2$.
\end{itemize}

System \eqref{model-equ} is  said to be at resonance in $l^2$ if
$\operatorname{dim}\ker (I-A\xi^{\alpha-1}) \geq1$, otherwise, it is said to be
non-resonant. The requirement $1\leq\operatorname{dim}\ker (I-A\xi^{\alpha-1})$
is to make the problem to be resonant and the requirement
$\operatorname{dim}\ker (I-A\xi^{\alpha-1}) <\infty$ is to make the
kernel operator to be a Fredholm operator which is a basic requirement in
applying the coincidence degree theory introduced by Mawhin.

In a recent paper \cite{Zhou2}, the authors studied the three-point BVPs \eqref{model-equ} at resonance in
infinite dimension space by assuming one of the following conditions holds

(A1)  $A\xi^{\alpha-1}$ is idempotent, that is,
$A^2\xi^{2\alpha-2}=A\xi^{\alpha-1}$;

(A2)  $A^2\xi^{2\alpha-2}=I$, where $I$ stands for the identity
operator in $l^2$.

The assumptions (A1) and (A2) are important in constructing the operator $Q$ in \cite{Zhou2}
which plays a key role in the process of the proof.
Our objective in this paper is
 to remove the restricted conditions (A1) and (A2) to study the existence of solutions for BVPs \eqref{model-equ}. It deserves to point out
that the problem is new even when $\alpha=2$, that is, system \eqref{model-equ} is a
second order differential system with resonant
boundary conditions.

In the past three decades, the existence of solutions
for the fractional differential equations with boundary value
conditions  have attained a great deal of attentions from many researchers, for instance,
 see \cite{Ba3,YCHen,JWH1,JWH2,Kosmatov1,Kosmatov2,Ba1,zhou1}.
However, all results derived in these papers are for one equation with
$\operatorname{dim}\ker L =1$ or  for two equations with
$\operatorname{dim}\ker L =2$. The case of problems where the $\operatorname{dim}\ker L$ can
take any value in $\mathbb{N}$  have treated with little attention.

 Recently, the authors in \cite{PhungTruong,Phung2}
investigated the following second order differential system
\begin{equation}
\begin{gathered}
u''(t)=f(t,u(t),u'(t)), \quad  0<t<1,\\
u'(0)=\theta,u(1)=Au(\eta)
\end{gathered}
\end{equation}
where $f:[0,1]\times\mathbb{R}^n\times\mathbb{R}^n\to\mathbb{R}^n$
is a Carath\'eodory function and the square matrix $A$
satisfies certain conditions.
These results for second order differential equations in
\cite{PhungTruong} and \cite{Phung2}  were generalized to fractional order case $\alpha\in(1,2]$ in \cite{GeZhouEJQTDE} and \cite{Zhou2}. It should be highlighted that, in \cite{Phung2}, the authors successfully removed  the constricted conditions used in \cite{PhungTruong} by making use of the property of Moore-Penrose pseudo-inverse matrix technique. Inspired by this, in this paper, we use the generalized inverse of the bounded linear operator in infinite dimensional space \cite{Petry} to remove the restricted conditions (A1) and (A2),
so that we can derive the existence of the solution for three-point BVPs \eqref{model-equ}.

We proceed as follows: in Section 2, we give some necessary
background and some preparations for our consideration.
The proof of the main results is presented
in Section 3 by applying the coincidence degree theory of
Mawhin.  In Section 4, an illustrative example is included.

\section{Preliminaries}

In this section, we present some necessary definitions and lemmas which
will be used later. These definitions and lemmas can be found in
 \cite{Mawhin1,Kilbas,Mawhin,Petry}  and the references therein.\\

\begin{definition}[\cite{Kilbas}] \label{def2.1}
The fractional integral of order $\alpha>0$ of a function
$x:(0,\infty)\to \mathbb{R}$ is given by
\[
I^{\alpha}_{0^+}x(t)=\frac{1}{\Gamma(\alpha)}
\int^{t}_{0}{(t-s)^{\alpha-1}x(s)ds},
\]
provided that the right-hand side is pointwise defined on $(0,\infty)$.
\end{definition}

\begin{remark}\label{rmk2.1}
The notation $I^{\alpha}_{0^+}x(t)|_{t=0}$ means that the limit is taken
at almost all points of the right-hand side neighborhood
$(0, \varepsilon) (\varepsilon > 0)$ of $0$ as follows:
\[
I^{\alpha}_{0^+}x(t)|_{t=0}=\lim_{t\to 0^+}I^{\alpha}_{0^+}x(t).
\]
Generally,  $I^{\alpha}_{0^+}x(t)|_{t=0}$ is not necessarily to be zero.
For instance, let $\alpha\in(0,1)$, $x(t)=t^{-\alpha}$.
Then
\[
I^{\alpha}_{0^+}t^{-\alpha}|_{t=0}=\lim_{t\to0^+}
\frac{1}{\Gamma(\alpha)}\int_0^t(t-s)^{\alpha-1}s^{-\alpha}ds
=\lim_{t\to0^+}\Gamma(1-\alpha)=\Gamma(1-\alpha).
\]
\end{remark}

\begin{definition}[\cite{Kilbas}]\label{def2.2}
  The fractional derivative of order $\alpha>0$
of a function $x:(0,\infty)\to \mathbb{R}$ in Riemann-Liouville sence is given by
\[
D^{\alpha}_{0^+}x(t)=\frac{1}{\Gamma(n-\alpha)}
\Big(\frac{d}{dt}\Big)^{n}\int^{t}_{0}{\frac{x(s)}{(t-s)^{\alpha-n+1}}ds},
\]
where $n=[\alpha]+1$, provided that the right-hand side is pointwise defined
on $(0,\infty)$.
\end{definition}

\begin{lemma}[\cite{Kilbas}] \label{Lem2.1}
Assume that $x\in C(0,+\infty)\cap L_{\rm loc}(0,+\infty)$ with a
fractional derivative of order $\alpha>0$ belonging to
$C(0,+\infty)\cap L_{\rm loc}(0,+\infty)$. Then
\[
I_{0^+}^{\alpha}D_{0^+}^{\alpha}x(t)=x(t)+c_1t^{\alpha-1}
+c_2t^{\alpha-2}+\dots+c_{n}t^{\alpha-n},
\]
for some $c_i\in\mathbb{R},i=1,\dots,n$, where $n=[\alpha]+1$.
\end{lemma}

For any $x(t)=(x_1(t),x_2(t),\dots)^\top \in l^2$, the
fractional derivative of order $\alpha>0$ of $x$ is defined by
\[
D^{\alpha}_{0^+}x(t)=(D^{\alpha}_{0^+}x_1(t),D^{\alpha}_{0^+}x_2(t),\dots)^\top
\in l^2.
\]
The following definitions and the coincidence degree theory are fundamental in the
proof of our main result. We refer the readers to \cite{Mawhin1,Mawhin}.\\

\begin{definition}\label{Fredhop}
 Let $X$ and $Y$ be normed spaces. A linear operator
$L:\operatorname{dom}(L)\subset X\to Y$ is said to be a Fredholm operator of
index zero provided that
\begin{itemize}
\item[(i)] $\operatorname{im}L$ is a closed subset of $Y$;

\item[(ii)] $\operatorname{dim}\ker L=\operatorname{codim}\operatorname{im}L<+\infty$.
\end{itemize}
\end{definition}

It follows from Definition \ref{Fredhop} that, if $L$ is a Fredholm operator of index zero, then there exist continuous projectors
 $P:X\to X$ and $Q:Y\to Y$ such that
\[
\operatorname{im}P=\ker L,\quad
\ker Q= \operatorname{im}L,\quad
X=\ker L\oplus\ker P,\quad Y=\operatorname{im}L\oplus\operatorname{im}Q,
\]
and the mapping $L|_{\operatorname{dom}L\cap\ker P}:\operatorname{dom}
L\cap\ker P\to\operatorname{im}L$ is invertible.
We denote the inverse of $L|_{\operatorname{dom}L\cap\ker P}$ by
$K_P:\operatorname{im}L\to\operatorname{dom}L\cap\ker P$.
The generalized inverse of $L$  denoted by
$K_{P,Q}:Y\to\operatorname{dom}L\cap\ker P$ is defined by
$K_{P,Q}=K_P(I-Q)$. Furthermore, for every isomorphism
$J:\operatorname{im}Q\to \ker L$, we can obtain that  the mapping
$K_{P,Q}+JQ: Y\to \operatorname{dom}L$ is also an isomorphism and for all
$x\in \operatorname{dom}L$, we have
\begin{equation}\label{isomorphism}
(K_{P,Q}+JQ)^{-1}x= (L+J^{-1}P)x.
\end{equation}

\begin{definition}\label{lcompact}
Let $L$ be a Fredholm operator of index zero,
$\Omega \subseteq X$ be a bounded subset and
$\operatorname{dom}L \cap \Omega \neq \emptyset$.
Then the operator $N: \overline{\Omega}\to Y$ is called to be
$L$-compact in $\overline{\Omega}$ if
\begin{itemize}
\item[(i)] the mapping $QN:\overline{\Omega} \to Y$ is continuous and
$QN(\overline{\Omega}) \subseteq Y $ is bounded;

\item[(ii)] the mapping $K_{P,Q}N:\overline{\Omega} \to X$ is
completely continuous.
\end{itemize}
\end{definition}

The following lemma is the main tool in this paper.
%

\begin{lemma}[\cite{Mawhin}] \label{lem2.2}
 Let $\Omega\subset X$ be a bounded subset, $L$ be a Fredholm mapping of index zero
and $N$ be $L$-compact in $\overline{\Omega}$. Suppose that the following
conditions are satisfied:
\begin{itemize}
\item[(i)] $Lx\neq\lambda Nx$ for every
$(x,\lambda)\in((\operatorname{dom}L\backslash\ker L)\cap\partial\Omega)\times(0,1)$;

\item[(ii)] $Nx\not\in \operatorname{im}L$ for every
$x\in\ker L\cap\partial\Omega$;

\item[(iii)] $\deg(JQN|_{\ker L\cap\partial\Omega},\Omega\cap\ker L,0)\neq0$,
with $Q:Y\to Y$ a continuous projector such that $\ker Q=\operatorname{im}L$ and
$J:\operatorname{im}Q\to\ker L$ is an isomorphism.
\end{itemize}
Then the equation $Lx=Nx$ has at least one solution in
$\operatorname{dom}L\cap\overline{\Omega}$.
\end{lemma}

In this paper, we use spaces $\mathbb{X}$, $\mathbb{Y}$ introduced as
\[
\mathbb{X}=\big\{x(t)\in l^2:x(t)=I_{0+}^{\alpha-1}u(t),
u\in C([0,1];l^2),t\in [0,1] \big\}
\]
with the norm
$\|x\|_\mathbb{X}=\max \{\|x\|_{C([0,1];l^2)},\|D_{0+}^{\alpha-1}x\|_{C([0,1];l^2)}\}$
and $\mathbb{Y}=L^1([0,1];l^2)$ with the norm
$\|y\|_{L^1([0,1];l^2)}=\int_0^1\|y(s)\|_{l^2}ds$,
respectively, where $\|x\|_{C([0,1];l^2)}=\sup_{t\in[0,1]}\|x(t)\|_{l^2}$.

\begin{lemma}\label{Lem2.3}
$F\subset \mathbb{X}$ is a sequentially compact set if and only if
$F(t)$ is relatively compact and equicontinuous which are understood in the
following sense:
\begin{itemize}
\item[(1)]  for any $t\in[0,1]$, $F(t):=\{x(t)|x\in F\}$ is  a relatively
compact set in $l^2$;
\item[(2)]  for any given $\varepsilon>0$, there
exists $\delta>0$ such that for all $x\in F$,
\[
\|x(t_1)-x(t_2)\|_{l^2}<\varepsilon,\;\|D_{0^+}^{\alpha-1}x(t_1)
-D_{0^+}^{\alpha-1}x(t_2)\|_{l^2}<\varepsilon,
\]
for $t_1,t_2\in[0,1]$, $|t_1-t_2|<\delta$.
\end{itemize}
\end{lemma}

In order to use Lemma \ref{lem2.2}, we define the linear operator
$L:\operatorname{dom}L \subseteq\mathbb{X}\to\mathbb{Y}$ by
\begin{equation}\label{Ldef}
 Lx:=D_{0^+ }^\alpha x,
\end{equation}
where $\operatorname{dom}L=\{x\in X:  D_{0^+ }^\alpha x\in Y,
I_{0+}^{2-\alpha}x(0)=\theta, x(1)=Ax(\xi)\}$ and  define $N:X\to Y$
by
\begin{equation}\label{N}
Nx(t):=f(t,x(t),D_{0^+}^{\alpha-1}x(t)), \quad  t\in [0,1].
\end{equation}
Then the problem \eqref{model-equ} can be equivalently rewritten as $Lx=Nx$.

Now we define operator $\mathcal{M}$ as:
\begin{equation}\label{M-def}
\mathcal{M}=I-A\xi^{\alpha-1},
\end{equation}
and define a continuous linear operator $h:\mathbb{Y}\rightarrow l^2$ by
\begin{equation}\label{hdef}
h(y):=\frac{A}{\Gamma(\alpha)}\int_0^{\xi}(\xi-s)^{\alpha-1}y(s)ds-\frac{I}{\Gamma(\alpha)}\int_0^1(1-s)^{\alpha-1}y(s)ds.
\end{equation}
In order to remove  the restricted conditions (A1) and (A2), we will employ the following lemma
on the property of bounded linear operator in general Hilbert space.
\begin{lemma}\label{M+exist}\cite{Petry}
If $\mathcal{M}$ is a bounded linear transformation from Hilbert space $H_1$ to Hilbert space $H_2$ with closed rang $R(\mathcal{M})$, then the generalized inverse $\mathcal{M}^+$ of $A$ is characterized as the unique solution $X$ of the following equivalent equations:\\
(\uppercase\expandafter{\romannumeral1}) $X\mathcal{M}X=X, \;(X\mathcal{M})^*=X\mathcal{M}, \;\mathcal{M}X\mathcal{M}=\mathcal{M}, (\mathcal{M}X)^*=\mathcal{M}X$;\\
(\uppercase\expandafter{\romannumeral2}) $\mathcal{M}X=P_{R(\mathcal{M})},\; N(X^*)=N(\mathcal{M})$ where $P_{R(\mathcal{M})}$ denotes the orthogonal projection of $H_2$ onto $R(\mathcal{M})$;\\
(\uppercase\expandafter{\romannumeral3}) $\mathcal{M}X=P_{R(\mathcal{M})},\; X\mathcal{M}=P_{R(\mathcal{M}^*)}, \; X\mathcal{M}X=X$;\\
(\uppercase\expandafter{\romannumeral4}) $X\mathcal{M}\mathcal{M}^*=\mathcal{M}^*,XX^*\mathcal{M}^*=X$;\\
(\uppercase\expandafter{\romannumeral5}) $X\mathcal{M}x=x$ for all $x\in R(\mathcal{M}^*)$ and $Xy=0$ for all $y\in N(\mathcal{M}^*)$;\\
(\uppercase\expandafter{\romannumeral6}) $X\mathcal{M}=P_{R(\mathcal{M}^*)}, \; N(X)=N(\mathcal{M}^*)$;\\
(\uppercase\expandafter{\romannumeral7}) $\mathcal{M}X=P_{R(\mathcal{M})},\;  X\mathcal{M}=P_{R(\mathcal{M})}$.
\end{lemma}

\begin{remark}
By the definition of $\mathcal{M}$ given in \eqref{M-def}, since $A$ is a bounded linear operator,
 we know that $\mathcal{M}$ satisfies the condition in Lemma \ref{M+exist}. Thus, for such $\mathcal{M}$, there exists
 unique $\mathcal{M}^+$ satisfies the equations in  Lemma \ref{M+exist}.
\end{remark}

The next lemma plays a vital role in estimating the boundedness of some sets.

\begin{lemma}\label{boundinequ}\cite{Zhou2}
Let $z_1,z_2\geq 0$, $\gamma_1,\gamma_2\in[0,1)$ and $\lambda_i,\mu_i\geq0,i=1,2,3$,
and the following two inequalities hold,
\begin{equation}\label{baseineq}
\begin{gathered}
z_1\leq\lambda_1z_1^{\gamma_1}+\lambda_2z_2+\lambda_3,\\
z_2\leq\mu_1z_1+\mu_2z_2^{\gamma_2}+\mu_3.
\end{gathered}
\end{equation}
Then $z_1,z_2$ is bounded if $\lambda_2\mu_1<1$.
\end{lemma}

\begin{lemma}\label{lemL}
The operator $L$, defined by \eqref{Ldef}, is a Fredholm operator of index zero.
\end{lemma}

\noindent
Proof\\
For any $x\in \operatorname{dom}L$,  by Lemma \ref{Lem2.1}
and  $I_{0+}^{2-\alpha}x(0)=\theta$, we have
\begin{equation}\label{integralequation}
x(t)=I_{0+}^\alpha Lx(t)+ct^{\alpha-1},
 \quad  c\in l^2,  \; t\in [0,1],
 \end{equation}
which, together with $x(1)=Ax(\xi)$, yields
\begin{equation}\label{KerLeq}
\begin{aligned}
\ker L&=\{x\in \mathbb{X}:x(t)=ct^{\alpha-1}, t\in [0,1], c\in \ker\\
&\simeq \ker\mathcal{M}t^{\alpha-1}.
\end{aligned}
\end{equation}
 Now we claim that
\begin{equation}\label{ImL}
\operatorname{im}L=\{y\in Y: h(y) \in \operatorname{im}\mathcal{M}\}.
\end{equation}
Actually, for any $y\in \operatorname{im}L$, there exists a function
$x\in \operatorname{dom}L$ such that $y=Lx$. It follows from
\eqref{integralequation} that
$x(t)=I_{0^+}^{\alpha}y(t)+ct^{\alpha-1}$, this jointly with $x(1)=Ax(\xi)$, follows
\[
\frac{A}{\Gamma(\alpha)}\int_0^\xi(\xi-s)^{\alpha-1}y(s)ds
-\frac{I}{\Gamma(\alpha)}\int_0^1(1-s)^{\alpha-1}y(s)ds
=\mathcal{M}c,\quad  c\in l^2,
\]
which means that  $h(y)\in \operatorname{im}\mathcal{M}$.

On the other hand, for any $y\in\mathbb{Y}$ satisfying
$h(y)\in \operatorname{im}\mathcal{M}$, there exists a constant
$c^*$ such that $h(y)=\mathcal{M}c^*$.
Let $x^*(t)=I_{0^+ }^\alpha y(t)+c^*t^{\alpha-1}$, a straightforward
computation shows that $x^*(0)=\theta$ and
$x^*(1)=Ax^*(\xi)$. Hence,
$x^*\in \operatorname{dom}L$ and $y(t)=D_{0^+ }^\alpha x^*(t)$, which implies
that $y\in \operatorname{im}L$.

 Furthermore, notice that if  $y=ct^{\alpha-1}$,
$c\in  l^2$, then
 \begin{equation}\label{gker}
\begin{aligned}
 h(y)
&=\frac{A}{\Gamma(\alpha)}\int_0^\xi(\xi-s)^{\alpha-1}y(s)ds
 -\frac{I}{\Gamma(\alpha)}\int_0^1(1-s)^{\alpha-1}y(s)ds\\
&=\frac{(\xi^{2\alpha-1}A-I)c}{\Gamma(\alpha)\Gamma(2\alpha)}.
\end{aligned}
 \end{equation}

Also, the relation
\begin{equation}\label{MM+}
(I-\mathcal{M}\mathcal{M}^+)(\xi^{2\alpha-1}A-I)=(\xi^{\alpha}-1)(I-\mathcal{M}\mathcal{M}^+)
\end{equation}
holds.
This is deduced from\\
$$(I-\mathcal{M}\mathcal{M}^+)(I-A\xi^{\alpha-1})=(I-\mathcal{M}\mathcal{M}^+)\mathcal{M}=0,$$
which is equivalent to
\begin{equation}\label{MM-equi}
\begin{aligned}
&(I-\mathcal{M}\mathcal{M}^+)A\xi^{\alpha-1}=(I-\mathcal{M}\mathcal{M}^+)\\
&\Leftrightarrow(I-\mathcal{M}\mathcal{M}^+)\xi^{2\alpha-1}A=\xi^{\alpha}(I-\mathcal{M}\mathcal{M}^+)\\
&\Leftrightarrow(I-\mathcal{M}\mathcal{M}^+)(\xi^{2\alpha-1}A-I)=(\xi^{\alpha}-1)(I-\mathcal{M}\mathcal{M}^+).
\end{aligned}
\end{equation}
Define the continuous linear mapping $Q:\mathbb{Y}\to\mathbb{Y}$ by
\begin{equation}\label{Q}
Qy(t):=\frac{\Gamma(\alpha)\Gamma(2\alpha)}{\xi^\alpha-1}(I-\mathcal{M}\mathcal{M}^+)h(y)t^{\alpha-1},
\quad t\in [0,1],\; \;y\in \mathbb{Y}.
\end{equation}
Then it follows from  \eqref{hdef}, \eqref{ImL}, \eqref{MM+} and Lemma \ref{M+exist} that
\begin{align*}
Q^2y(t)
&=\frac{\Gamma(\alpha)\Gamma(2\alpha)}{\xi^\alpha-1}(I-\mathcal{M}\mathcal{M}^+)h(Qy(t))t^{\alpha-1}\\
&=\frac{\Gamma(\alpha)\Gamma(2\alpha)}{\xi^\alpha-1}(I-\mathcal{M}\mathcal{M}^+)
 \frac{(\xi^{2\alpha-1}A-I)}{\Gamma(\alpha)\Gamma(2\alpha)}
\frac{\Gamma(\alpha)\Gamma(2\alpha)}{\xi^\alpha-1}(I-\mathcal{M}\mathcal{M}^+)h(y)t^{\alpha-1}\\
&=\frac{\Gamma(\alpha)\Gamma(2\alpha)}{\xi^\alpha-1}(I-\mathcal{M}\mathcal{M}^+)^2h(y)t^{\alpha-1}\\
&=\frac{\Gamma(\alpha)\Gamma(2\alpha)}{\xi^\alpha-1}(I-\mathcal{M}\mathcal{M}^+)h(y)t^{\alpha-1}
=Qy(t),
\end{align*}
and
\begin{gather*}
y\in \ker  Q\Leftrightarrow h(y)\in\ker  (I-\mathcal{M}\mathcal{M}^+)
\Leftrightarrow h(y)\in \operatorname{im} \mathcal{M}\mathcal{M}^+
\Leftrightarrow h(y)\in \operatorname{im} \mathcal{M}
\Leftrightarrow y\in \operatorname{im}L,
\end{gather*}
which implies that $Q$ is a projection operator with $\ker Q= \operatorname{im}L$.
Therefore,
$\mathbb{Y}=\ker Q \oplus \operatorname{im} Q=\operatorname{im}L \oplus
\operatorname{im}Q $.\\
 Finally, we shall prove that $\operatorname{im}Q=\ker L$. Indeed, for any
$z\in \operatorname{im} Q$, let $z=Qy$, $y\in \mathbb{Y}$.
we have
\[
(\mathcal{M}\mathcal{M}^+)z(t)=(\mathcal{M}\mathcal{M}^+)Qy(t)
=\frac{\Gamma(\alpha)\Gamma(2\alpha)}{\xi^\alpha-1}\mathcal{M}\mathcal{M}^+(I-\mathcal{M}\mathcal{M}^+)g(y)t^{\alpha-1}
=\theta,
\]
which implies $z\in \ker  L$.
Conversely, for each $z\in \ker L$, there exists a constant
$c^*\in \ker (\mathcal{M})$ such that $z=c^*t^{\alpha-1}$ for
$t\in [0,1]$. By \eqref{gker} and \eqref{MM-equi}, we derive
\[
Qz(t)=\frac{\Gamma(\alpha)\Gamma(2\alpha)}{\xi^\alpha-1}
(I-\mathcal{M}\mathcal{M}^+)h(c^*t^{\alpha-1})t^{\alpha-1}=c^*t^{\alpha-1}=z(t),\quad  t\in [0,1],
\]
which implies that $z\in \operatorname{im}Q$. Hence we know that
$\operatorname{im}Q=\ker L$. By assumption that $\operatorname{dim}\ker (I-A\xi^{\alpha-1})<\infty$,
the operator $L$ is a Fredholm operator of
index zero. The proof is completed.

Now to establish the generalized inverse of $L$, we define the operator $P:\mathbb{X}\to\mathbb{X}$ by
\begin{equation}\label{Popdef}
Px(t)=\frac{1}{\Gamma(\alpha)}(I-\mathcal{M}^+\mathcal{M})D_{0+}^{\alpha-1}x(0)t^{\alpha-1}.
\end{equation}

\begin{lemma}\label{lem-Pdef}
The following assertions hold:\\
(1) The mapping $P:\mathbb{X}\to\mathbb{X}$ defined by \eqref{Popdef},
is a continuous projector satisfying
\[
\operatorname{im}P=\ker  L, \quad \mathbb{X}=\ker  L\oplus\ker  P;
\]
(2) The linear operator $K_P: \operatorname{im}L\to \operatorname{dom}L \cap\ker  P$, which is the inverse of $L|_{\operatorname{dom}L \cap\ker  P}$,
can be written as
\begin{equation}\label{KPydef}
K_Py(t)=\mathcal{M}^+h(y)t^{\alpha-1}+I_{0+}^\alpha y(t),
\end{equation}
moreover, $K_P$ satisfies
\[
\|K_Py\|_{\mathbb{X}}\leq C\|y\|_{L^1([0,1];l^2)},
\]
where $C=1+\|\mathcal{M}^+\mathcal{M}\|(1+\|A\|).$
\end{lemma}

\noindent
Proof\\
(1) By Lemma \ref{M+exist}, $I-\mathcal{M}^+\mathcal{M}$ is a projection on $\ker \mathcal{M}\subset l^2$.
It follows from \eqref{Popdef} that $P$ is a continuous projection.
If $v\in\operatorname{im}P$, there exists  $x\in \mathbb{X}$ such that $v=Px$,
then
\[
v=\frac{1}{\Gamma(\alpha)}
(I-\mathcal{M}^+\mathcal{M})D_0^{\alpha-1}x(0)t^{\alpha-1}.
\]
By \eqref{KerLeq} and  Lemma \ref{M+exist},
 we have $\mathcal{M}(I-\mathcal{M}^+\mathcal{M})D_0^{\alpha-1}x(0)=0$,
which gives $v\in \ker  L$.
Conversely, if $v\in \ker L$, then
$v(t)=c_*t^{\alpha-1}$ for some $c_*\in\ker  \mathcal{M}=\operatorname{im}(I-\mathcal{M}^+\mathcal{M})$, that is,
$c_*=(I-\mathcal{M}^+\mathcal{M})\bar{c}_*$ for $\bar{c}_*\in l^2$. Thus, we deduce that
\begin{align*}
Pv(t)&=\frac{1}{\Gamma(\alpha)}(I-\mathcal{M}^+\mathcal{M})D_{0+}^{\alpha-1}v(0)t^{\alpha-1}
=(I-\mathcal{M}^+\mathcal{M})c_*t^{\alpha-1}\\
&=(I-\mathcal{M}^+\mathcal{M})^2\bar{c}_*t^{\alpha-1}=(I-\mathcal{M}^+\mathcal{M})\bar{c}_*t^{\alpha-1}\\
&=v(t), \quad t\in [0,1],
\end{align*}
which gives $v\in\operatorname{im}P$. Thus, we get that
 $\ker L=\operatorname{im}P$ and consequently
$\mathbb{X}={\rm ker }L\oplus \ker P$.

(2) Let $y\in \operatorname{im}L$. There exists
$x\in \operatorname{dom}L$ such that $y=Lx$ and $h(y)\in \operatorname{im}\mathcal{M}$. By the
definitions of $P$ and $K_P$, we obtain that
\begin{align*}
PK_Py(t)&=\frac{1}{\Gamma(\alpha)}(I-\mathcal{M}^+\mathcal{M})D_{0^+}^{\alpha-1}(K_Py(0))t^{\alpha-1}\\
&=\frac{1}{\Gamma(\alpha)}(I-\mathcal{M}^+\mathcal{M})[\mathcal{M}^+D_{0^+}^{\alpha-1}h(y(0))+I_{0^+}^\alpha y(0)]t^{\alpha-1}\\
&=0,
\end{align*}
and
\[
\mathcal{M}(K_Py(0))=\mathcal{M}[\mathcal{M}^+h(y(0))+I_{0^+}^{\alpha}y(0)]=h(y).
\]
Thus, $K_Py\in {\rm ker }P \cap\operatorname{dom}L$, $K_P$ is well defined.

On the other hand, if $x\in {\rm ker }P \cap\operatorname{dom}L$, then $x(t)=I_{0^+}^{\alpha}Lx(t)+ct^{\alpha-1}$, and
$$\mathcal{M}c=h(Lx),\ \ c\in {\rm ker }(I-\mathcal{M}^+\mathcal{M}).$$
Hence
\begin{align*}
K_PL_Px(t)&=\mathcal{M^+}h(Lx)t^{\alpha-1}+I_{0^+}^{\alpha}Lx(t)\\
&=\mathcal{M}^+h(Lx)\\
&=\mathcal{M}^+\mathcal{M}ct^{\alpha-1}+I_{0^+}^{\alpha}Lx(t)\\
&=ct^{\alpha-1}+I_{0^+}^{\alpha}Lx(t)=x(t)
\end{align*}
and $L_PK_Px(t)=x(t), \ t\in [0,1]$ for all $x\in \operatorname{im}L$, then $K_P=L_P^{-1}$.

Finally, by the definition of $K_P$, we have
\begin{equation}\label{DKPy}
(D^{\alpha-1}_{0^+}K_Py)(t)=\Gamma(\alpha)\mathcal{M}^+h(y)+I_{0^+}^1y(t).
\end{equation}
It follows from \eqref{hdef}, \eqref{KPydef} and \eqref{DKPy} that
\begin{align*}
\|D_{0^+}^{\alpha-1}K_Py\|_{C([0,1];l^2)}&= \Gamma(\alpha)\|\mathcal{M}^+\|\|h(y)\|_{C([0,1];l^2)}+\bigg\|\int_0^{\cdot}y(s)ds\bigg\|_{C([0,1];l^2)}\\
&\leq \Gamma(\alpha)\|\mathcal{M}^+\|\|h(y)\|_{C([0,1];l^2)}+\|y\|_{L^1([0,1];l^2)},
\end{align*}
\begin{align*}
\|K_Py\|_{C([0,1];l^2)}&= \|\mathcal{M}^+\|\|h(y)\|_{C([0,1];l^2)}
+\bigg\|\frac{1}{\Gamma(\alpha)}\int_0^{\cdot}(\cdot-s)^{\alpha-1}y(s)ds\bigg\|_{C([0,1];l^2)}\\
&\leq \|\mathcal{M}^+\|\|h(y)\|_{C([0,1];l^2)}+\frac{1}{\Gamma(\alpha)}\|y\|_{L^1([0,1];l^2)},
\end{align*}
and
\begin{align*}
\|h(y)\|_{C([0,1];l^2)}\leq \frac{1}{\Gamma(\alpha)}(1+\|A\|)\|y\|_{L^1([0,1];l^2)}.
\end{align*}
This show that
\begin{align*}
\|K_Py\|_\mathbb{X}&=\max\{\|D_{0^+}^{\alpha-1}K_Py\|_{C([0,1];l^2)},\|K_Py\|_{C([0,1];l^2)}\}\\
&\leq[1+\|\mathcal{M}^+\mathcal{M}\|(1+\|A\|)]\|y\|_{L^1}.
\end{align*}
This completes of the proof.

\begin{lemma} \label{lem-N-compact}
Let $f$ be a Carath\'eodory function. Then  $N$, defined by \eqref{N},
is L-compact.
\end{lemma}

\noindent
Proof\\
Let $\Omega$ be a bounded subset in $\mathbb{X}$. By
hypothesis (iii) on the function $f$,
 there exists a function $\varphi_\Omega(t)\in L^1[0,1]$ such that for all
$x\in \Omega$,
\begin{equation}\label{fphi}
\|f(t,x(t),D_{0^+}^{\alpha-1}x(t))\|_{l^2}\leq \varphi_\Omega(t), \quad
\text{a.e. }t\in [0,1],
\end{equation}
which, along with  \eqref{hdef}, implies
\begin{equation}\label{gbound}
\begin{aligned}
 \|h(Nx(t))\|_{l^2}
&=\bigg\|\frac{A}{\Gamma(\alpha)}\int_0^\xi(\xi-s)^{\alpha-1}
 {f(s,x(s),D_{0^+}^{\alpha-1}x(s))}ds\\
&\quad -\frac{I}{\Gamma(\alpha)}\int_0^1(1-s)^{\alpha-1}
 {f(s,x(s),D_{0^+}^{\alpha-1}x(s))}ds\bigg\|_{l^2}\\
&\leq \frac{\|A\|+1}{\Gamma(\alpha)}\|\varphi_\Omega\|_{L^1[0,1]}.
\end{aligned}
\end{equation}
Thus,  from \eqref{Q} and \eqref{gbound}, it follows that
\begin{equation}
\begin{aligned}
  \|QNx\|_{L^1([0,1];l^2)}
&=\big\|\frac{\Gamma(\alpha)\Gamma(2\alpha)}{\xi^\alpha-1}(I-\mathcal{M}^+\mathcal{M})h(Nx)
 \big\|_{l^2}\int_0^1s^{\alpha-1}ds\\
&\leq  \frac{\Gamma(2\alpha)(\|A\|+1)\|I-\mathcal{M}^+\mathcal{M}\|}{|1-\xi^\alpha|}
 \|\varphi_\Omega\|_{L^1[0,1]}<\infty.
 \end{aligned}
 \end{equation}
This shows that $QN(\overline{\Omega}) \subseteq\mathbb{Y}$ is bounded. The
continuity of $QN$ follows from the hypothesis on $f$ and the
Lebesgue dominated convergence theorem.

 Next, we shall show that $K_{P,Q}N$ is completely continuous.
For any $x\in \Omega$, we have
\begin{equation}\label{KPQdef}
\begin{aligned}
K_{P,Q}Nx(t)
&=K_P(I-Q)Nx(t)=K_PNx(t)-K_PQNx(t)\\
&=\mathcal{M}^+h(Nx)t^{\alpha-1}+I_{0+}^\alpha Nx(t)\\
&\;\;\;-
\frac{\Gamma(\alpha)\Gamma(2\alpha)}{\xi^\alpha-1}(I-\mathcal{M}^+\mathcal{M})h(Nx(t))I_{0+}^\alpha
t^{\alpha-1},
\end{aligned}
\end{equation}
and
\begin{equation}\label{KPQdefD}
\begin{aligned}
 D_{0^+}^{\alpha-1}K_{P,Q}Nx(t)&=\Gamma(\alpha)\mathcal{M}^+h(Nx)+I_{0+}^1Nx(t)\\
 &\;\;\;-
\frac{\Gamma(\alpha)\Gamma(2\alpha)}{\xi^\alpha-1}(I-\mathcal{M}^+\mathcal{M})h(Nx(t))I_{0+}^1
t^{\alpha-1}.
\end{aligned}
\end{equation}
By the hypothesis on $f$ and the Lebesgue dominated convergence theorem,
it is easy to see that $K_{P,Q}N$ is continuous.
Since $f$ is a Carath\'eodory function, for every bounded set
$\Omega_0\subseteq l^2\times l^2$,
the set $\{f(t,u,v):(u,v)\in \Omega_0\}$ is relatively compact set in $l^2$.
 Therefore, for almost all  $t\in[0,1]$, $\{K_{P,Q}Nx(t):x\in\Omega\}$ and
 $\{D_{0^+}^{\alpha-1}K_{P,Q}Nx(t):x\in\Omega\}$ are relatively compact in $l^2$.

From \eqref{gbound}, \eqref{KPQdef} and \eqref{KPQdefD}, we derive that
 \begin{align*}
&\ \ \ \ \|K_{P,Q}Nx\|_{C([0,1];l^2)}\\
&=\big\|\mathcal{M}^+h(Nx)t^{\alpha-1}+I_{0+}^\alpha Nx(t)-\frac{\Gamma(\alpha)\Gamma(2\alpha)}{\xi^\alpha-1}
(I-\mathcal{M}^+\mathcal{M})h(Nx(t))I_{0+}^\alpha  t^{\alpha-1}\big\|_{C([0,1];l^2)}\\
&\leq\|\mathcal{M}^+\|\frac{\|A\|+1}{\alpha\Gamma(\alpha)}\|\varphi_{\Omega}\|_{L^1[0,1]}+\frac{1}{\Gamma(\alpha)}\|\varphi_\Omega\|_{L^1(0,1)}
+\frac{\Gamma(2\alpha)\|I-\mathcal{M}^+\mathcal{M}\|}{|\xi^\alpha-1|}\|h(Nx(t))\|_{l^2}\\
&\leq \|\mathcal{M}^+\|\frac{\|A\|+1}{\alpha\Gamma(\alpha)}\|\varphi_{\Omega}\|_{L^1[0,1]}+\frac{1}{\Gamma(\alpha)}\|\varphi_\Omega\|_{L^1(0,1)}
 +\frac{\Gamma(2\alpha)\|I-\mathcal{M}^+\mathcal{M}\|(\|A\|+1)}{\Gamma(\alpha)|\xi^\alpha-1|}
\|\varphi_\Omega\|_{L^1(0,1)}\\
&<\infty,
\end{align*}
and
\begin{align*}
&\ \ \ \ \|D_{0^+}^{\alpha-1}K_{P,Q}Nx\|_{C([0,1];l^2)}\\
&\leq\frac{\Gamma(\alpha)}{\alpha}\|\mathcal{M}^+\| \|\varphi_{\Omega}\|_{L^1}+\big\|I_{0+}^1 Nx(t)-\frac{\Gamma(\alpha)\Gamma(2\alpha)}{\xi^\alpha-1}
 (I-\mathcal{M}^+\mathcal{M})h(Nx(t))I_{0+}^1  t^{\alpha-1}\big\|_{C([0,1];l^2)}\\
&\leq \frac{\Gamma(\alpha)}{\alpha}\|\mathcal{M}^+\| \|\varphi_{\Omega}\|_{L^1}+\|\varphi_\Omega\|_{L^1(0,1)}+\frac{\Gamma(2\alpha)
\|I-\mathcal{M}^+\mathcal{M}\|}{|\xi^\alpha-1|}\|h(Nx(t))\|_{l^2}\\
&\leq \frac{\Gamma(\alpha)}{\alpha}\|\mathcal{M}^+\| \|\varphi_{\Omega}\|_{L^1}+\|\varphi_\Omega\|_{L^1(0,1)}
 +\frac{\Gamma(2\alpha)\|I-\mathcal{M}^+\mathcal{M}\|(\|A\|+1)}{\Gamma(\alpha)|\xi^\alpha-1|}
\|\varphi_\Omega\|_{L^1(0,1)}\\
&<\infty,
\end{align*}
which shows that $K_{P,Q}N\overline{\Omega}$ is uniformly bounded in
$\mathbb{X}$.  Noting  that
\begin{equation}\label{inequality}
 b^p-a^p\leq (b-a)^p \quad \text{for any } b\geq a>0, 0< p\leq1.
\end{equation}
 for any $t_1,t_2\in [0,1]$ with $t_1<t_2$,  we shall see that
\begin{align*}
&\ \ \ \ \|K_{P,Q}Nx(t_2)-K_{P,Q}Nx(t_1)\|_{l^2}\\
&\leq\frac{1}{\Gamma(\alpha)}\Big\|\Gamma(\alpha) \mathcal{M}^+h(Nx)(t_2^{\alpha-1}-t_1^{\alpha-1})+\int_0^{t_1}{[(t_2-s)^{\alpha-1}-(t_1-s)^{\alpha-1}]
 Nx(s)}ds\\
&\quad +\int_{t_1}^{t_2}{(t_2-s)^{\alpha-1}Nx(s)}ds
-\frac{\Gamma(\alpha)\Gamma(2\alpha)}{\xi^\alpha-1}(I-\mathcal{M}^+\mathcal{M})h(Nx(t))
 [I_{0+}^\alpha t_2^{\alpha-1}-I_{0+}^\alpha t_1^{\alpha-1}]\Big\|_{l^2}\\
&\leq\|\mathcal{M}^+h(Nx)\|_{l^2}(t_2-t_1)^{\alpha-1}+ \frac{1}{\Gamma(\alpha)}\int_0^{t_1}{(t_2-t_1)^{\alpha-1}\varphi_\Omega(s)}ds
+\frac{1}{\Gamma(\alpha)}\int_{t_1}^{t_2}{\varphi_\Omega(s)}ds\\
&\quad +\frac{\Gamma^2(\alpha)\|I-\mathcal{M}^+\mathcal{M}\|(\|A\|+1)}{|\xi^\alpha-1|}
\|\varphi_\Omega\|_{L^1(0,1)}|t_2^{2\alpha-1}- t_1^{2\alpha-1}|\\
&\leq \frac{1}{\Gamma(\alpha)}\int_0^{t_1}(t_2-t_1)^{\alpha-1}\varphi_{\Omega}(s)ds+\frac{1}{\Gamma(\alpha)}\int_{t_1}^{t^2}\varphi_{\Omega}(s)ds\\
&\quad+\|\mathcal{M}^+\|\frac{\|A\|+1}{\Gamma(\alpha)}\|\varphi_{\Omega}\|_{L^1}(t_2-t_1)^{\alpha-1}\\
&\quad+\frac{\Gamma^2(\alpha)\|I-\mathcal{M}^+\mathcal{M}\|(\|A\|+1)}{|\xi^{\alpha}-1|}\|\varphi_{\Omega}\|_{L^1}|t_2^{2\alpha-1}-t_1^{2\alpha-1}|\to 0, \quad \text{as } t_2\to t_1
\end{align*}
and
\begin{align*}
&\ \ \ \ \|D_{0^+}^{\alpha-1}K_{P,Q}Nx(t_2)-D_{0^+}^{\alpha-1}K_{P,Q}Nx(t_1)\|_{l^2}\\
&=\big\|\int_{t_1}^{t_2}{Nx(s)}ds\big\|_{l^2}
+\big\|\frac{\Gamma(\alpha)\Gamma(2\alpha)}{\xi^\alpha-1}(I-\mathcal{M}^+\mathcal{M})h(Nx(t))
 \int_{t_1}^{t_2}s^{\alpha-1}ds\big\|_{l^2}\\
&\leq \int_{t_1}^{t_2}{\varphi_\Omega(s)}ds
+\frac{\Gamma(2\alpha)\|I-\rho_A\|(\|A\|+1)}{|\xi^\alpha-1|}
\|\varphi_\Omega\|_{L^1(0,1)}|t_2^{\alpha}-t_1^{\alpha}|
\to 0, \quad \text{as }  t_2\to t_1.
\end{align*}
Then  $K_{P,Q}N\overline{\Omega}$ is equicontinuous in
$\mathbb{X}$. By Lemma \ref{Lem2.3}, $K_{P,Q}N\overline{\Omega}\subseteq\mathbb{X}$
is  relatively compact. Thus we can conclude that the operator $N$
is $L$-compact in $\overline{\Omega}$.  The proof is completed.

\section{Main results}

\begin{theorem}\label{mainresults}
Let $f$ be a Carath\'eodory function and the following conditions
hold:
\begin{itemize}
\item[(H1)] There exist five nonnegative functions $a_1, a_2, b_1, b_2,c\in L^1[0,1]$
and two constants $\gamma_1,\gamma_2\in[0,1)$ such that
for all $t\in [0,1]$, $u,v \in l^2$,
\[
\|f(t,u,v)\|_{l^2}\leq a_1(t)\|u\|_{l^2}+b_1(t)\|v\|_{l^2}
+a_2(t)\|u\|^{\gamma_1}_{l^2}+b_2(t)\|v\|^{\gamma_2}_{l^2} +c(t)
\]
holds.
\item[(H2)] There exists a constant $A_1>0$ such that for
$x\in\operatorname{dom}L$, if $\|D_{0^+}^{\alpha-1}x(t)\|_{l^2}>A_1$
for all $t\in [0,1]$, then
\begin{align*}
& \frac{A}{\Gamma(\alpha)}\int_0^\xi(\xi-s)^{\alpha-1}{f(s,x(s),
 D_{0^+}^{\alpha-1}x(s))}ds\\
&-\frac{I}{\Gamma(\alpha)}\int_0^1(1-s)^{\alpha-1}{f(s,x(s),
 D_{0^+}^{\alpha-1}x(s))}ds\notin\operatorname{im} \mathcal{M}.
\end{align*}

\item[(H3)]
 There exists a  constant $A_2>0$ and an isomorphism $J:\operatorname{im}Q\to \ker L$ such that for any $e=\{(e_i)\}\in l^2$
satisfying $e=\xi^{\alpha-1}Ae$ and $\|e\|_{l^2}>A_2$, either
\[
\langle e,JQNe\rangle_{l^2}\leq 0 \quad \text{or} \quad
\langle e,JQNe\rangle_{l^2}\geq 0\quad \text{holds},
\]
where $\langle \cdot,\cdot\rangle_{l^2}$ is the inner product in $l^2$.
\end{itemize}
Then  \eqref{model-equ} has at least one solution in space $\mathbb{X}$ provided that
\begin{equation}\label{Thcondition}
\begin{gathered}
\Gamma(\alpha)>\max\big\{(\|I-\mathcal{M}^+\mathcal{M}\|+1)\|a_1\|_{L^1(0,1)},
(\|I-\mathcal{M}^+\mathcal{M}\|+1)\|b_1\|_{L^1(0,1)}\big\},\\
\frac{(\|I-\mathcal{M}^+\mathcal{M}\|+1)^2\|a_1\|_{L^1(0,1)}\|b_1\|_{L^1(0,1)}}
{(\Gamma(\alpha)-(\|I-\mathcal{M}^+\mathcal{M}\|+1)\|a_1\|_{L^1(0,1)})
(\Gamma(\alpha)-(\|I-\mathcal{M}^+\mathcal{M}\|+1)\|b_1\|_{L^1(0,1)})}<1.
\end{gathered}
\end{equation}
\end{theorem}

To prove the above theorem, we need the following auxiliary lemmas.
\begin{lemma}\label{lem-Ome1}
The set
$
\Omega_1=\big\{x\in \operatorname{dom}L\backslash\ker
L:Lx=\lambda Nx \text{ for some } \lambda \in [0,1]\}
$ is bounded in $\mathbb{X}$.
\end{lemma}

\noindent
Proof\\
For any $x\in \Omega_1$, $x \notin\ker L$, we have
$\lambda\neq 0$. Since $Nx\in \operatorname{im}L=\ker Q$,
 by \eqref{ImL}, we have $h(Nx)\in \operatorname{im} \mathcal{M}$, where
\begin{equation}\label{hDef}
\begin{aligned}
h(Nx)&=\frac{A}{\Gamma(\alpha)}\int_0^\xi(\xi-s)^{\alpha-1}{f(s,x(s),
 D_{0^+}^{\alpha-1}x(s))}ds\\
&\quad -\frac{I}{\Gamma(\alpha)}\int_0^1(1-s)^{\alpha-1}{f(s,x(s),
 D_{0^+}^{\alpha-1}x(s))}ds.
\end{aligned}
\end{equation}
From (H2), there exists $t_0\in [0,1]$ such that
 $\|D_{0^+}^{\alpha-1}x(t_0)\|_{l^2}\leq A_1$. Then from the equality
$D_{0+}^{\alpha-1}x(0)=D_{0+}^{\alpha-1}x(t_0)-\int_0^t{D_{0+}^\alpha
x(s)}ds$, we deduce that
\[
\|D_{0+}^{\alpha-1}x(0)\|_{l^2}
\leq A_1+\|D_{0+}^\alpha x\|_{L^1(0,1;l^2)}
= A_1+\|Lx\|_1\leq A_1+\|Nx\|_{L^1(0,1;l^2)},
\]
which implies
\begin{equation}\label{P}
\|Px\|_{\mathbb{X}}
=\|\frac{1}{\Gamma(\alpha)}(I-\mathcal{M}^+\mathcal{M})D_{0+}^{\alpha-1}x(0)t^{\alpha-1}
\|_{\mathbb{X}}
\leq \frac{\|I-\mathcal{M}^+\mathcal{M}\|}{\Gamma(\alpha)}( A_1+\|Nx\|_{L^1(0,1;l^2)}).
\end{equation}
Further,  for $x\in \Omega_1$, since $\operatorname{im}P=\ker L,
\mathbb{X}=\ker  L\oplus\ker P $, we have
$(I-P)x \in \operatorname{dom}L \cap \ker P$ and $LPx=\theta$. Then
\begin{equation}\label{I-P}
\begin{aligned}
\|(I-P)x\|_{\mathbb{X}}
&=\|K_PL(I-P)x\|_{\mathbb{X}}\leq \|K_PLx\|_{\mathbb{X}}\\
&\leq \frac{1}{\Gamma(\alpha)}\|Lx\|_{L^1(0,1;l^2)}
 \leq \frac{1}{\Gamma(\alpha)}\|Nx\|_{L^1(0,1;l^2)}.
\end{aligned}
\end{equation}
From \eqref{P} and \eqref{I-P}, we conclude that
\begin{equation}
\begin{aligned}\label{|x|}
\|x\|_{\mathbb{X}}
&=\|Px+(I-P)x\|_{\mathbb{X}} \leq \|Px\|_{\mathbb{X}}+\|(I-P)x\|_{\mathbb{X}}\\
&\leq\frac{\|I-\mathcal{M}^+\mathcal{M}\|}{\Gamma(\alpha)}A_1+\frac{\|I-\mathcal{M}^+\mathcal{M}\|+1}{\Gamma(\alpha)}
 \|Nx\|_{L^1(0,1;l^2)}.
\end{aligned}
\end{equation}
Moreover, by the definition of $N$ and (H1), we derive
\begin{equation}\label{NX}
 \begin{aligned}
&\ \ \ \ \|Nx\|_{L^1(0,1;l^2)}\\
&= \int_0^1{\|f(s,x(s),D_{0^+}^{\alpha-1}x(s))\|_{l^2}}dt\\
&\leq\|a_1\|_{L^1(0,1)}\|x\|_{C([0,1];l^2)}
 +\|b_1\|_{L^1(0,1)}\|D_{0^+}^{\alpha-1}x\|_{C([0,1];l^2)}\\
&\quad +\|a_2\|_{L^1(0,1)}\|x\|_{C([0,1];l^2)}^{\gamma_1}
 +\|b_2\|_{L^1(0,1)}\|D_{0^+}^{\alpha-1}x\|_{C([0,1];l^2)}^{\gamma_2}
+\|c\|_{L^1(0,1)}.
\end{aligned}
\end{equation}
Thus,
\begin{equation}\label{x}
 \begin{aligned}
\|x\|_{\mathbb{X}}
&\leq \frac{\|I-\mathcal{M}^+\mathcal{M}\|}{\Gamma(\alpha)}A_1+\frac{\|I-\mathcal{M}^+\mathcal{M}\|+1}{\Gamma(\alpha)}
\Big(\|a_1\|_{L^1(0,1)}\|x\|_{C([0,1];l^2)}\\
&\quad +\|b_1\|_{L^1(0,1)}\|D_{0^+}^{\alpha-1}x\|_{C([0,1];l^2)}\Big)
+\frac{\|I-\mathcal{M}^+\mathcal{M}\|+1}{\Gamma(\alpha)} \\
&\ \ \quad\times \Big(\|a_2\|_{L^1(0,1)}\|x\|_{C([0,1];l^2)}^{\gamma_1}
 +\|b_2\|_{L^1(0,1)}\|D_{0^+}^{\alpha-1}x\|_{C([0,1];l^2)}^{\gamma_2}\\
&\ \ \ \ \quad +\|c\|_{L^1(0,1)}\Big).
\end{aligned}
\end{equation}
It follows from \eqref{Thcondition}, \eqref{x},
$\|x\|_{C([0,1];l^2)}\leq \|x\|_{\mathbb{X}}$,
$\|D_{0^+}^{\alpha-1}x\|_{C([0,1];l^2)}\leq \|x\|_{\mathbb{X}}$ and Lemma \ref{boundinequ} that there exists $M_0>0$ such that
\[
\max\{\|x\|_{C([0,1];l^2)}, \|D_{0^+}^{\alpha-1}x\|_{C([0,1];l^2)}\}\leq M_0,
\]
which means that $\Omega_1$ is bounded in $\mathbb{X}$.

\begin{lemma}\label{lem-Ome2}
The set
$
\Omega_2=\{x \in\ker L:N x\in  \operatorname{im}L\}
$ is bounded in $\mathbb{X}$.
\end{lemma}

\noindent
Proof\\
For any $x\in \Omega_2$, it follows from $x \in\ker L $ that
$x=et^{\alpha-1}$ for some $e\in \ker \mathcal{M}\subset l^2$,
and it follows from
$N x\in  \operatorname{im}L$ that $h(N x) \in \operatorname{im}\mathcal{M}$,
where $h(N x)$ is defined by \eqref{hDef}. By hypothesis $(H_2)$, we arrive at
$\|D_{0^+}^{\alpha-1}x(t_0)\|_{l^2}=\|e\|_{l^2}\Gamma(\alpha)\leq A_1$.
Thus, we obtain $\|x\|\leq \|e\|_{l^2}\Gamma(\alpha)\leq A_1$.
i.e.,  $\Omega_2$ is bounded in $\mathbb{X}$.

\begin{lemma}\label{lem-Ome3}
Let $\Omega_3=\{x\in\ker  L:-\lambda x+(1-\lambda)JQNx=\theta, \; \lambda \in [0,1]\}$ if the first part of $(H_3)$ holds,
and $\Omega_3=\{x\in \ker L:\lambda x+(1-\lambda)JQNx=\theta, \; \lambda \in [0,1]\}$ if the other part of (H3) holds.
Then, the set $\Omega_3$ is bounded in $\mathbb{X}$.
\end{lemma}

\noindent
Proof\\
If the first part of $(H_3)$ holds, that is, $\langle e, JQNe\rangle_{l^2}\leq 0$,

then for any $x\in \Omega_3$, we know that
\[
x=et^{\alpha-1}\quad \text{with } e\in \ker \mathcal{M}\text{ and }
\lambda x=(1-\lambda)JQNx.
\]
If $\lambda=0$, we have $Nx\in\ker Q= \operatorname{im}L$, then
$x\in \Omega_2$, by the argument above, we get that $\|x\| \leq A_1$.
Moreover, if $\lambda\in (0,1]$ and if $\|e\|_{l^2}>A_2$, by
(H3), we deduce that
\[
0<\lambda \|e\|^2_{l^2}=\lambda\langle e,e\rangle_{l^2}
=(1-\lambda)\langle e,JQNe\rangle_{l^2}\leq 0,
\]
which is a contradiction. Then
$\|x\|_{\mathbb{X}}=\|et^{\alpha-1}\|_{\mathbb{X}}\leq
\max\{\|e\|_{l^2}, \Gamma(\alpha)\|e\|_{l^2}\}$.
That is to say, $\Omega_3$  is bounded.\\
For the case of the second part of $(H_3)$ holds, we can obtain the result that $\Omega_3$ is bounded by a similar method as above, so we omit it.

{\it Proof of Theorem 3.1}:
We first construct an open bounded subset $\Omega$ in $\mathbb{X}$ such that
$\cup_{i = 1}^3\overline{\Omega}_i \subseteq \Omega$.
By Lemmas \ref{lemL} and \ref{lem-N-compact}, we know that $L$ is a Fredholm
operator of index zero and $N$ is $L$-compact on $\overline{\Omega}$. Thus,
it follows from Lemmas \ref{lem-Ome1}, \ref{lem-Ome2} and \ref{lem-Ome3} that
conditions (i) and (ii) of Lemma \ref{lem2.2} hold.
By the construction of $\Omega$ and the argument
above,  to complete the theorem, it suffices to prove that
condition (iii) of Lemma \ref{lem2.2} is satisfied.
To this end, let
\begin{equation}
H(x,\lambda)=\pm \lambda x+(1-\lambda)JQNx,
\end{equation}
here we let the isomorphism $J:\operatorname{im}Q\to \ker L$ be the
identical operator. Since $\Omega_3\subseteq \Omega$,
$H(x,\lambda)\neq 0$ for $(x,\lambda)\in\ker L \cap
\partial \Omega \times [0,1]$, then by homotopy property of degree, we obtain
\begin{align*}
\deg \left( {JQN{|_{\ker L\cap \partial \Omega }},\Omega\cap\ker L,\theta} \right)
&= \deg \left( {H\left(
\cdot,0\right),\Omega\cap\ker L,\theta} \right)\\
&=\deg \left( {H\left( \cdot,1\right),\Omega\cap\ker L,\theta} \right)\\
&=\deg \left( {\pm Id,\Omega\cap\ker L,\theta} \right)=\pm 1\neq 0.
\end{align*}
Thus (H3) of Lemma \ref{lem2.2} is fulfilled and Theorem
\ref{mainresults} is proved.

\section{Example}

In this section, we shall present an example to illustrate our main result
in $l^2$.

Consider the following system  with  $\operatorname{dim}\ker L=k$,
$k=1,2,3,\dots$ in $l^2$.
\begin{equation}\label{exampleker2}
\begin{gathered}
\begin{aligned}
&D_{0^+}^{3/2} \begin{pmatrix}
x_1(t)\\
x_2(t)\\
x_3(t)\\
x_4(t)\\
\vdots
\end{pmatrix} \\
&=\frac{1}{10} \begin{pmatrix}
\begin{cases}
 1, &\text{if } \|D_{0^+}^{1/2}x(t)\|_{l^2}< 1\\
 D_{0^+}^{1/2}x_1(t)+[D_{0^+}^{1/2}x_1(t)]^{-1}-1,
 &\text{if } \|D_{0^+}^{1/2}x(t)\|_{l^2}\geq 1
\end{cases} \\
\big(x_2(t)+D_{0^+}^{1/2}x_2(t)\big)/2\\
\big(x_3(t)+D_{0^+}^{1/2}x_3(t)\big)/2^2\\
\big(x_4(t)+D_{0^+}^{1/2}x_4(t)\big)/2^3\\
\vdots
\end{pmatrix}
\end{aligned}\\
 I_{0^+}^{\frac{1}{2}}x_i(0)=0, \quad i=1,2,\dots\\
x(1)=Ax(1/4).
\end{gathered}
\end{equation}
For all $t\in[0, 1]$, let
$u=(x_1,x_2,x_3,\dots)^{\top}$, $v=(y_1,y_2,y_3,\dots)^{\top}\in l^2$ and
$f=(f_1,f_2,\dots)^{\top}$ with
\[
f_1(t,u,v)=\begin{cases}
 1/10, &\text{if } \|v\|_{l^2}< 1,\\
 (y_1+y_1^{-1}-1)/10, &\text{if } \|v\|_{l^2}\geq 1,
\end{cases}
\]
$f_i(t,u,v)=\frac{1}{5}\frac{x_i+y_i}{2^{i}}$, $i=2,3,4,\dots$.
Moreover,
 \begin{equation}\label{Adef2x2x}
A=\begin{bmatrix}
{B_1}&0&{0}&{0}&{0}&{0}&\dots\\  {0}&{B_2}&0&{0}&{0}&{0}&\dots\\
\vdots& &\ddots& & & &\vdots\\
{0}&0&{0}&{B_k}&{0}&{0}&{\dots}\\
 {0}&0&{0}&{0}&{0}&{0}&{\dots}\\{0}&0&{0}&{0}&{0}&{0}&{\dots}
 \\ \vdots& & & & & &\ddots
 \end{bmatrix}
\quad\text{with}\quad
 B_i=\begin{bmatrix} {\frac{3}{2}}&{0}&{0}\\ {0}&{\frac{7}{4}}&{0}\\{0}&{0}&{2}
\end{bmatrix},
\end{equation}
and we denote
\[\mathcal{M}_i=I-\xi^{\frac{1}{2}}B_i=
\begin{pmatrix}
\frac{1}{4}&0&0\\
0&\frac{1}{8}&0\\
0&0&0
\end{pmatrix},
\]
then
\[
\mathcal{M}^+_i=
\begin{pmatrix}
4&0&0\\
0&8&0\\
0&0&0
\end{pmatrix},
\]
$i=1,2,\dots,k$, $k\in \mathbf{N}$.
Obviously, $\operatorname{dim}\ker (I_3-\xi^{\alpha-1}B_i)=\operatorname{dim}\ker (I_3-B_i/2)=1$,
 $i=1,2,\dots$, where $I_3$ is the  $3\times 3$ identity matrix.
Then $\operatorname{dim}\ker (I-A\xi^{\alpha-1})=k$,
$k\in \mathbf{N}$ and the problem \eqref{exampleker2},
with $A$ and $f$ defined above,  has one solution if and only if problem
\eqref{model-equ} admits one solution.

Checking (H1) of Theorem \ref{mainresults}:
For some $r\in \mathbb{R}$,
 $\Omega=\{(u,v)\in l^2\times l^2:
\|u\|_{l^2}\leq r, \|v\|_{l^2}\leq r\}$, let
$\varphi_\Omega(t)=\frac{1}{10}[(r+1/r+1)^2+\frac{4r^2}{3}]^{1/2}\in L^1[0,1]$.
Letting
\begin{equation}\label{abex2e}
  a_1(t)=b_1(t)=\frac{1}{5\sqrt{3}},\quad
 a_2(t)=b_2(t)=0 , \quad c(t)=\frac{r+1/r+1}{10},
\end{equation}
 condition (H1) is satisfied.

Checking (H2) of Theorem \ref{mainresults}:
From the definition of $f$ it follows that $f_1>1/10>0$ when
$ \|D_{0^+}^{1/2}x(t)\|_{l^2}> 1$.
This,
\[
(B_1\xi^{\alpha}-I)\begin{pmatrix}
{f_1}\\
{f_2}\\
{f_3}
\end{pmatrix}
=\begin{bmatrix}
{-13/16}&{0}&{0}\\
{0}&{-25/32}&{0}\\
{0}&{0}&{-3/4}
\end{bmatrix}
\begin{pmatrix}
{f_1}\\
{f_2}\\
{f_3}
\end{pmatrix}=\begin{pmatrix}
-\frac{13}{16}f_1\\
-\frac{25}{32}f_2\\
-\frac{3}{4}f_3
\end{pmatrix},
\]
and $\operatorname{im} (\mathcal{M})
=\{(2\tau_1,\tau_1,0,2\tau_2,\tau_2,0,2\tau_3,\tau_3,0\dots)^{\top}:\tau_{i}\in \mathbb{R},
i=1,2,\dots\}$ implies that
 condition (H2) is satisfied.

Checking (H3) of Theorem \ref{mainresults}:
Since $\operatorname{dim}\ker (\mathcal{M})=k$, $k\in \mathbf{N}$,
for any $e\in l^2$  satisfying $e=\xi^{\alpha-1}Ae$, $e$ can be expressed as
$e=e_1+e_2+\dots+e_k$,  with
\[
e_i=\sigma_{i}(\varepsilon_{3i}), \quad \sigma_{i} \in \mathbb{R},\quad
 i=1,2,\dots,k,j=1,2,
\]
where $\varepsilon_j=( {0,0, \dots 0,\mathop 1_{{\rm{j}} - th}, 0, 0,\dots }
 )^{\top}\in l^2$ is a vector with all elements equaling to 0 except the
$j$-th equaling to 1, $j=1,2,\dots$, that is
\[
e=(0,0,\sigma_{1},0,0,\sigma_{2},0,0,\sigma_3,\cdots)^{\top}.
\]

In addition, for any $y\in \mathbb{Y}$, by \eqref{Q}, we have
\begin{equation}\label{Qexample11}
Qy(t)=\frac{\Gamma(\alpha)\Gamma(2\alpha)}{\xi^\alpha-1}
(I-\mathcal{M}\mathcal{M}^+)h(y)t^{\alpha-1}=-\frac{8\sqrt{\pi}}{7}(I-\mathcal{M}\mathcal{M}^+)h(y)t^{\alpha-1},
\end{equation}
where
\begin{equation}
h(y)= \frac{A}{\Gamma(\alpha)}\int_0^{1/4}(\frac{1}{4}-s)^{\alpha-1}y(s)ds
-\frac{I}{\Gamma(\alpha)}\int_0^1(1-s)^{\alpha-1}y(s)ds.
\end{equation}
By \eqref{N}, let $d=t^{1/2}+\frac{\sqrt{\pi}}{2}$, we have
\begin{equation}\label{Nct1ex11}
N(et^{1/2})\\
=\frac{1}{10} \begin{cases}
\Big(1,0, \frac{d\sigma_{1}}{2^2},0,0,\frac{d\sigma_{2}}{2^{5}},0,0,\frac{d\sigma_3}{2^8},\dots,0,0,\frac{d\sigma_i}{2^{3i-1}}\dots\Big)^\top,\\
\quad\text{if } |\sigma_{11}|< 1,\; 2\leq i\leq k;\\[4pt]
 \Big(-1, 0, \frac{d\sigma_{1}}{2^2},0,0,\frac{d\sigma_{2}}{2^{5}},0,0,\frac{d\sigma_3}{2^8},\dots,0,0,\frac{d\sigma_i}{2^{3i-1}}\dots\Big)^\top_,\\
\quad\text{if }|\sigma_{11}|\geq 1,\; 2\leq i\leq k.
 \end{cases}
\end{equation}
For $|\sigma_{11}|>1$, let $\hat{d}=\frac{\pi}{128}+\frac{\sqrt{\pi}}{24}$, $\tilde{d}=\frac{\pi}{8}+\frac{\sqrt{\pi}}{3}$, and let $A_2=1$, we have
\begin{align*}
\int_0^{1/4}(1/4-s)^{1/2}(Nes^{1/2})ds=\frac{1}{10}\Big(
-\frac{1}{12},0,\frac{\hat{d}\sigma_1}{2^2},0,0,\frac{\hat{d}\sigma_2}{2^5},0,0,\frac{\hat{d}\sigma_3}{2^8},0,0,\frac{\hat{d}\sigma_4}{2^{11}},
\cdots\Big)^{\top},
\end{align*}
and
\begin{align*}
\int_0^{1}(1-s)^{1/2}(Nes^{1/2})ds=\frac{1}{10}\Big(
-\frac{3}{2},0,\frac{\tilde{d}\sigma_1}{2^2},0,0,\frac{\tilde{d}\sigma_2}{2^5},0,0,\frac{\tilde{d}\sigma_3}{2^8},0,0,\frac{\tilde{d}\sigma_4}{2^{11}},
\cdots\Big)^{\top},
\end{align*}

\begin{align*}
h(Net^{1/2})&=\frac{A}{\Gamma(\alpha)}\int_0^{1/4}(1/4-s)^{1/2}Nes^{1/2}ds-\frac{I}{\Gamma(\alpha)}\int_0^1(1-s)^{1/2}Nes^{1/2}ds\\
=\frac{1}{5\sqrt{\pi}}&\Big(
\frac{11}{8},0,\frac{(2\hat{d}-\tilde{d})\sigma_1}{2^2},0,0,\frac{(2\hat{d}-\tilde{d})\sigma_2}{2^5},0,0,\frac{(2\hat{d}-\tilde{d})\sigma_3}{2^8},
0,0,\frac{(2\hat{d}-\tilde{d})\sigma_4}{2^{11}},\cdots
\Big)^{\top}.
\end{align*}
Then
\begin{align*}
Q(Net^{1/2})
=&-\frac{8\sqrt{\pi}t^{1/2}}{7}(I-\mathcal{M}\mathcal{M}^+)h(Net^{1/2})\\
=&-\frac{8t^{1/2}}{35}\Big(
0,0,\frac{(2\hat{d}-\tilde{d})\sigma_1}{2^2},0,0,\frac{(2\hat{d}-\tilde{d})\sigma_2}{2^5},0,0,\frac{(2\hat{d}-\tilde{d})\sigma_3}{2^8},
0,0,\frac{(2\hat{d}-\tilde{d})\sigma_4}{2^{11}},
\cdots\Big)^\top,
\end{align*}
and
\begin{align*}
\langle e,QNet^{1/2}\rangle
=&-\frac{8t^{1/2}}{35}
\Big(\frac{(2\hat{d}-\tilde{d})}{2^2}\sigma_1^2+\frac{(2\hat{d}-\tilde{d})}{2^5}\sigma_2^2+\frac{(2\hat{d}-\tilde{d})}{2^8}\sigma_3^2
+\frac{(2\hat{d}-\tilde{d})}{2^{11}}\sigma_4^2+\cdots
\Big)>0.\\
\end{align*}

Therefore, \eqref{exampleker2} admits at least one solution.

\begin{remark}
By a simply calculation of $B$, we can get
 \[
 B^2=
 \begin{pmatrix}
 \frac{9}{4}&0&0\\
 0&\frac{49}{16}&0\\
 0&0&4
\end{pmatrix},
 \]
from which, we can see that $A$ does not satisfies the conditions (A1) and (A2), so the result in \cite{Zhou2} is no longer applicable. Thus, our result is more general than the  one in  \cite{Zhou2}.
 \end{remark}

\section{Concluding remarks}
In this paper, we consider the fractional BVPs at resonance in $l^2$. The dimension of the kernel of the fractional differential operator with the boundary conditions be any positive integer.
We remove  the restricted conditions $A^2\xi^{2\alpha-2}=A\xi^{\alpha-1}$ and  $A^2\xi^{2\alpha-2}=I$
on the operator $A$, which have been used in \cite{Zhou2}.
Our result can also be easily generalized to other fractional BVPs, for instance,
\begin{equation}\label{consys}
\begin{gathered}
D_{0^+}^{\alpha}x(t)=f(t,x(t),D_{0^+}^{\alpha-1}x(t)), {\kern 5pt}1<\alpha\leq 2,{\kern 5pt} t\in (0,1),\\
x(0)=\theta, D_{0^+ }^{\alpha -1}x(1)=AD_{0^+}^{\alpha-1}x(\xi),
\end{gathered}
\end{equation}
where the bounded linear operator $A\in \mathcal{L}(l^2)$ satisfies $1\leq\operatorname{dim}\ker (I-A)<\infty$ which leads this system is resonant. Moreover, notice that $\mathbb{R}^n$ is the closed space of $l^2$,
taking $\alpha=2$, the system \eqref{consys} becomes the system of second order differential equations, which can be regarded as a generalization  results in \cite{PhungTruong} and \cite{Phung2}.

\section*{Acknowledgements} This work is supported by "the Fundamental Research Funds for the Central Universities" (No. CUSF-DH-D-2017083).

\end{document}